\documentclass[12pt]{article}
\usepackage{amssymb, amsmath}
\begin{document}
\renewcommand{\theequation}{\arabic{section}.\arabic{equation}}
\newtheorem{theorem}{Theorem}[section]
\newtheorem{lemma}{Lemma}[section]
\newtheorem{pro}{Proposition}[section]
\newtheorem{cor}{Corollary}[section]
\newcommand{\n}{\nonumber}
\newcommand{\tv}{\tilde{v}}
\newcommand{\tw}{\tilde{\omega}}
\renewcommand{\t}{\theta}
\newcommand{\w}{\omega}
\newcommand{\e}{\varepsilon}
\renewcommand{\a}{\alpha}
\renewcommand{\l}{\lambda}
\newcommand{\vare}{\varepsilon}
\newcommand{\s}{\sigma}
\renewcommand{\o}{\omega}
\newcommand{\bb}{\begin{equation}}
\newcommand{\ee}{\end{equation}}
\newcommand{\bq}{\begin{eqnarray}}
\newcommand{\eq}{\end{eqnarray}}
\newcommand{\bqn}{\begin{eqnarray*}}
\newcommand{\eqn}{\end{eqnarray*}}
\title{Nonexistence of  self-similar singularities for
the 3D incompressible Euler equations}
\author{Dongho Chae\thanks{The work was supported
partially by the KOSEF Grant no. R01-2005-000-10077-0.} \\
Department of Mathematics\\
              Sungkyunkwan University\\
              Suwon 440-746, Korea\\
  e-mail: {\it chae@skku.edu}}
 \date{}
\maketitle
\begin{abstract}
We prove that there exists no self-similar finite time blowing up
solution to the 3D incompressible Euler equations.  By similar
method we also show nonexistence of self-similar blowing up
solutions to the divergence-free transport equation in $\Bbb R^n$.
This result has  direct applications to the density dependent Euler
equations, the Boussinesq system, and the quasi-geostrophic
equations, for which we also show nonexistence of
 self-similar blowing up solutions.
\end{abstract}

\section{Incompressible Euler equations}
 \setcounter{equation}{0}
We are  concerned here on the following Euler equations for the
homogeneous incompressible fluid flows in $\Bbb R^3$.
 \[
\mathrm{ (E)}
 \left\{ \aligned
 &\frac{\partial v}{\partial t} +(v\cdot \nabla )v =-\nabla p ,
 \quad (x,t)\in {\Bbb R^3}\times (0, \infty) \\
 &\quad \textrm{div }\, v =0 , \quad (x,t)\in {\Bbb R^3}\times (0,
 \infty)\\
  &v(x,0)=v_0 (x), \quad x\in \Bbb R^3
  \endaligned
  \right.
  \]
where $v=(v_1, v_2, v_3 )$, $v_j =v_j (x, t)$, $j=1,2,3$, is the
velocity of the flow, $p=p(x,t)$ is the scalar pressure, and $v_0
$ is the given initial velocity, satisfying div $v_0 =0$.
 There are
well-known results on the local existence of  classical
solutions(see e.g. \cite{ maj1, kat1, con2} and references
therein). The problem of finite time blow-up of the local
classical solution is one of the most challenging open problem in
mathematical fluid mechanics. On this direction there is a
celebrated result on the blow-up criterion by Beale, Kato and
Majda(\cite{bea}). By geometric type of consideration some of the
possible scenarios to the possible singularity has been
excluded(see \cite{con4, cor3, den}. One of the main purposes of
this paper is to exclude the possibility of self-similar type of
singularities for the Euler system.\\
The  system (E) has scaling property that
  if $(v, p)$ is a
solution of the system (E), then for any $\lambda >0$ and $\alpha
\in \Bbb R $ the functions
 \bb
 \label{self}
  v^{\lambda, \alpha}(x,t)=\lambda ^\alpha v (\lambda x, \l^{\a +1}
  t),\quad p^{\l, \a}(x,t)=\l^{2\a}p(\l x, \l^{\a+1} t )
  \ee
  are also solutions of (E) with the initial data
  $ v^{\lambda, \alpha}_0(x)=\lambda ^\alpha v_0
   (\lambda x)$.
 In view of the scaling
  properties in (\ref{self}), the  self-similar blowing up
  solution $v(x,t)$ of (E) should be of the form,
  \bq
  \label{vel}
 v(x, t)&=&\frac{1}{(T_*-t)^{\frac{\a}{\a+1}}}
V\left(\frac{x}{(T_*-t)^{\frac{1}{\a+1}}}\right)
 \eq
 for  $\a \neq -1$ and $t$  sufficiently
 close to $T_*$.
  Substituting (\ref{vel}) into (E), we find that $V$ should be
  a solution of the system
\[
 (SE)
 \left\{
 \begin{aligned}
 \frac{\a}{\a +1} V+&\frac{1}{\a+1} (x\cdot\nabla ) V
 +(V\cdot \nabla )V =-\nabla  P,
 \\
  \textrm{div }\, V& =0
  \end{aligned}
  \right.
  \]
  for some scalar function $P$,
  which could be regarded as the Euler version of the Leray
  equations introduced in \cite{ler}. The question of existence of
  nontrivial solution to (SE) is equivalent to the  that of existence of
  nontrivial self-similar finite time blowing up
  solution to the Euler system of the form (\ref{vel}).
  Similar question for the
  3D Navier-Stokes equations was raised by J. Leray in \cite{ler}, and
  answered  negatively by the authors of \cite{nec}, the result of which was
  refined later in \cite{tsa}.
 Combining the energy conservation with a simple scaling argument,
  the author of this article showed that if there exists a nontrivial self-similar finite time
blowing up solution, then its helicity should be zero(\cite{cha1},
see also \cite{pom} for other related discussion).
  To the author's knowledge, however, the  possibility of self-similar blow-up of the
  form (\ref{vel})
  has never been excluded previously.
   In particular, due to lack of the laplacian term
  in the right hand side of the first equations of (SE), we cannot
  apply the argument of the maximum principle, which was crucial in the work \cite{nec}
   for the 3D Navier-Stokes equations.
  Using a completely different argument from those  used in \cite{cha1}, or \cite{nec},
  we prove here that there cannot be  self-similar blowing up solution to (E) of the form
  (\ref{vel}), if the vorticity decays sufficiently fast near infinity.
\begin{theorem}
There exists no finite time blowing up self-similar solution to
the 3D Euler equations of the form
  (\ref{vel}) for $t \in (0 ,T_*)$ with $\a \neq -1$, if there
  exists $p_1 >0$ such that
  the vorticity $\Omega=$curl $V\in L^p
 (\Bbb R^3;
 \Bbb R^3)$ for all $p\in (0, p_1)$.
\end{theorem}

\noindent{\it Remark 1.1} For example, if $\Omega\in
L^{1}_{loc}(\Bbb R^3;\Bbb R^3)$ and  there exist constants $R, K$
and $\delta >0$ such that $|\Omega(x)|\leq K e^{-\delta |x|}$ for
$|x|>R$, then we have $\Omega\in L^p (\Bbb R^3;\Bbb R^3)$ for all
$p\in (0, 1)$. Indeed, for all $p\in (0, 1)$, we have
 \bqn
\int_{\Bbb R^3} |\Omega(x)|^p dx
&=& \int_{|x|\leq R} |\Omega(x)|^p dx +\int_{|x|>R}  |\Omega (x)|^p \,dx\\
&\leq&|B_R |^{1-p}\left(\int_{|x|\leq R} |\Omega(x)| dx\right)^{p} +
K^p \int_{\Bbb R^3} e^{-p\delta |x|}dx <\infty ,
 \eqn
where $|B_R|$ is the volume of the ball $B_R$ of radius $R$.\\
\ \\
\noindent{\it Remark 1.2.} We note that there is no integrability
condition imposed on the velocity $V$ itself in the above theorem.
In particular, $V$  does not need to decay at infinity. For example,
if curl $V=\Omega $ has compact support in $\Bbb R^3$ with
$\mathrm{div}\, V=0$, we have by the Biot-Savart law,
$$V(x) =\frac{1}{4\pi} \int_{\Bbb R^3}
\frac{(x-y)\times \Omega (y)}{|x-y|^3}dy +\nabla h(x)
$$
for a harmonic function $h(x)$ in $\Bbb R^3$.
 Choosing
 $h(x)=2x_1^2-x_2^2-x_3^2$, for example, we have  $V(x)$, which grows
 to infinite in the positive $x_1$-direction.\\

The proof of Theorem 1.1 will follow  as a corollary of the
following more general theorem.

\begin{theorem}
 Suppose there exists $T>0$ such that we have  a
 representation of the vorticity of the solution, $v\in C([0,T);C^1 (\Bbb R^3 ;\Bbb R^3))$,
  to the 3D Euler
 equations by
 \bb\label{thm12}
 \o (x,t) = \Psi(t)\Omega (\Phi(t)x) \qquad \forall t\in [0, T)
 \ee
 where $\Psi(\cdot)\in C([0, T );(0, \infty))$,
 $\Phi(\cdot)\in C([0, T );\Bbb R^{3\times 3})$
 with $\mathrm{det}(\Phi(t))\neq 0$ on $[0, T)$;  $\Omega = \mathrm{curl}\, V$ for some
 $V$, and   there exists $p_1 >0$ such that
 $\Omega \in L^p
 (\Bbb R^3;
 \Bbb R^3)$ for all $p\in (0, p_1 )$.  Then,
necessarily
 either $\mathrm{det}(\Phi(t))\equiv \mathrm{det}(\Phi(0))$ on $[0, T)$, or
 $\Omega=0$.
\end{theorem}
{\bf Proof.} By consistency with the initial condition, $ \o_0
(x)=\Psi (0)\Omega (\Phi (0)x)$, and hence
  $\Omega (x)=\Psi(0)^{-1} \o_0 ([\Phi (0) ]^{-1} x)$ for all
$x\in \Bbb R^3$. We can rewrite the representation (\ref{thm12})
in the form,
 \bb \label{thm12a}
 \o (x,t) = G(t)\o_0 (F(t)x)
 \qquad \forall t\in [0, T),
 \ee
where $G(t)=\Psi (t)/\Psi (0)$, $F(t)=[\Phi (0)]^{-1}\Phi (t)$. In
order to prove the theorem it suffices to show that either
det$(F(t))=1$ for all $t\in [0,T)$, or $\o_0 =0$, since
det$(F(t))$= det$(\Phi (t))/$det$(\Phi(0))$.
   Let $a\mapsto X(a,t)$ be the
particle
 trajectory mapping, defined by the ordinary differential equations,
 $$
 \frac{\partial X(a,t)}{\partial t} =v(X(a,t),t) \quad;\quad X(a,0)=a.
 $$
We set $A(x,t):=X^{-1} (x ,t)$, which is called the back to label
map(\cite{con1}), satisfying
 \bb\label{inverse}
 A (X(a,t),t)=a, \quad X(A
 (x,t),t)=x.
 \ee
Taking curl of the first equation of (E), we obtain the vorticity
evolution equation,
$$
 \frac{\partial \o }{\partial t} +(v\cdot \nabla )\o =(\o
\cdot
 \nabla )v .
 $$
This, taking dot product with $\o$, leads to
 \bb\label{new0}
 \frac{\partial |\o |}{\partial t} +(v\cdot \nabla )|\o |=\a |\o|,
 \ee
where $\a (x,t)$ is defined as
 $$
 \a (x,t)=\left\{ \aligned &\sum_{i,j=1}^3 S_{ij}
 (x,t)\xi_i(x,t)\xi_j (x,t) &\quad \mbox{if} \quad\o (x,t)\neq 0\\
 &\qquad 0 &\quad \mbox{if} \quad \o (x,t)=0\endaligned \right.
 $$
 with
 $$ S_{ij} =\frac12 \left( \frac{\partial v_j}{\partial x_i}
 +\frac{\partial v_i}{\partial x_j}\right), \quad\mbox{and}\quad
\xi (x,t)=\frac{\o(x,t)}{|\o (x,t)|}.
$$
We note that (\ref{new0}) was previously derived in \cite{con4}. In
terms of the particle trajectory mapping we can rewrite (\ref{new0})
as
 \bb\label{new1}
 \frac{\partial }{\partial t} |\o
(X(a,t),t)|=\alpha (X(a,t),t) |\o (X(a,t),t)|.
 \ee
Integrating (\ref{new1}) along the particle trajectories $\{
X(a,t)\}$, we have \bb\label{new2}
 |\o (X(a,t),t)|=|\o_0 (a)|\exp \left[ \int_0 ^t \a
(X(a,s),s) ds \right].
  \ee
Taking into account the simple estimates
$$ -\|\nabla v(\cdot ,t)\|_{L^\infty}\leq \a
(x,t) \leq \|\nabla v (\cdot ,t)\|_{L^\infty} \quad \forall x\in
\Bbb R^3,
 $$
we obtain from (\ref{new2}) that
 \bqn
 \lefteqn{|\o_0 (a)|\exp \left[- \int_0 ^t \|\nabla v
      (\cdot,s)\|_{L^\infty} ds \right]\leq |\o (X(a,t),t)|}\hspace{1.in}\\
      && \qquad \leq |\o_0 (a)|\exp \left[ \int_0 ^t \|\nabla v
      (\cdot,s)\|_{L^\infty} ds \right],
 \eqn
 which, using the back to label map, can be rewritten as
  \bq\label{new3}
 \lefteqn{|\o_0 (A(x,t))|\exp \left[- \int_0 ^t \|\nabla v
      (\cdot,s)\|_{L^\infty} ds \right]\leq |\o
      (x,t)|}\hspace{1.in}\n\\
     && \leq |\o_0 (A(x,t))|\exp \left[ \int_0 ^t \|\nabla v
      (\cdot,s)\|_{L^\infty} ds \right].
 \eq
 Combining this with the self-similar representation formula in (\ref{thm12a}), we have
 \bq
  \label{new4}
 \lefteqn{|\o_0 (A(x,t))|\exp \left[- \int_0 ^t \|\nabla v
      (\cdot,s)\|_{L^\infty} ds \right]\leq G(t) |\o _0 (F(t)x)|}\hspace{1.3in}\n \\
      &&\leq
      |\o_0 (A(x,t))|\exp \left[ \int_0 ^t \|\nabla v
      (\cdot,s)\|_{L^\infty} ds \right].
 \eq
 Given $p\in (0, p_1)$, computing $L^p(\Bbb R^3)$ norm of the each side of (\ref{new4}),
 we derive
 \bq\label{new4a}
  \lefteqn{\|\o_0 \|_{L^p} \exp \left[- \int_0 ^t \|\nabla v
      (\cdot,s)\|_{L^\infty} ds \right]\leq G(t)[\mathrm{det}(F(t))]^{-\frac1p}
      \|\o _0 \|_{L^p} }\hspace{1.5in}\n \\
      &&\leq
      \|\o_0 \|_{L^p}\exp \left[ \int_0 ^t \|\nabla v
      (\cdot,s)\|_{L^\infty} ds \right],
\eq
 where we used the fact $\mathrm{det} (\nabla A (x,t))\equiv 1$.
 Now, suppose $\Omega \neq 0$, which is equivalent to assuming that $\o_0 \neq
 0$, then we divide (\ref{new4a}) by $\|\o_0 \|_{L^p}$ to obtain
 \bq \label{new5}
 \lefteqn{\exp \left[- \int_0 ^t \|\nabla v
      (\cdot,s)\|_{L^\infty} ds \right]\leq G(t)[\mathrm{det}(F(t))]^{-\frac1p}
      }\hspace{1.5in}\n \\
   &&\leq \exp \left[ \int_0 ^t \|\nabla v
      (\cdot,s)\|_{L^\infty} ds \right].
 \eq
 If there exists $t_1\in (0, T)$ such that
 $\mathrm{det}(F(t_1))\neq1$, then either $\mathrm{det}(F(t_1))>1$
 or $\mathrm{det}(F(t_1))<1$. In either case, setting $t=t_1$  and
  passing $p\searrow 0$ in
 (\ref{new5}), we deduce that
 $$\int_0 ^{t_1}\|\nabla v
      (\cdot,s)\|_{L^\infty} ds =\infty.
      $$
      This contradicts with the assumption that the flow is smooth on $(0, T)$, i.e $v\in C
      ([0,T); C^1 (\Bbb R^3;\Bbb R^3 ))$. $\square$\\
\ \\
\noindent{\bf Proof of Theorem 1.1} We apply Theorem 1.2 with
$$\Phi (t)=(T_* -t)^{-\frac{1}{\a +1}}I,\quad\mbox{and}\quad
 \Psi (t) =(T_* -t) ^{-1},
 $$
 where $I$ is the unit matrix in $\Bbb R^{3\times 3}$.
If $\a\neq -1$ and $t\neq 0$, then
 $$ \mathrm{det}(\Phi (t))=(T_* -t)^{-\frac{3}{\a +1}}
 \neq T_* ^{-\frac{3}{\a +1}} =\mathrm{det}(\Phi (0)).
 $$
 Hence, we conclude  that $\Omega=0$
by  Theorem 1.2.  In this case, there is no finite time blow-up
for $v(x,t)$,
 since the vorticity is zero. $\square$

\section{Divergence-free transport equation}
 \setcounter{equation}{0}

The previous argument in the proof of Theorem 1.1 can also be
applied to the following transport equations by a divergence-free
vector field in $\Bbb R^n$, $n\geq 2$.
 $$
  (TE)\left\{ \aligned
  &\frac{\partial \theta}{\partial t} +(v\cdot \nabla )\theta =0,\\
 &\quad \mathrm{ div}\, v =0, \\
  &\theta (x,0)=\theta_0 (x),
  \endaligned
  \right.
  $$
 where $v=(v_1, \cdots ,v_n )=v(x,t)$, and $\theta =\theta (x,t)$.
In view of the invariance of the transport equation under
 the scaling transform,
 \bqn
  &&v(x,t)\mapsto v^{\l, \a }(x,t)=\lambda ^\alpha v(\lambda x, \l^{\a +1}
  t),\\
  &&\theta(x,t)\mapsto \t ^{\l, \a,\beta} (x,t)= \l ^{\beta} \t (\lambda x, \l^{\a +1}
  t)
 \eqn
  for all $\a, \beta\in \Bbb R$ and $\lambda >0$, the self-similar
  blowing up solution is of the form,
 \bq\label{veqa}
 v(x, t)&=&\frac{1}{\left(T_*-t\right)^{\frac{\a}{\a+1}}}
V\left(\frac{x}{(T_*-t)^{\frac{1}{\a+1}}}\right),\\
\label{teqa}
 \theta  (x,t)&=&\frac{1}{\left(T_*-t\right)^{\beta}}
 \Theta \left(\frac{x}{(T_*-t)^{\frac{1}{\a+1}}}\right)
 \eq
 for $\a \neq -1$ and $t$ sufficiently close to $T_*$.
Substituting (\ref{veqa}) and (\ref{teqa}) into the above transport
equation, we obtain
 $$(ST)\left\{
 \aligned
 &\beta \Theta  +\frac{1}{\a +1} (x\cdot \nabla ) \Theta +(V\cdot
 \nabla )\Theta =0,\\
 &\qquad\mathrm{div}\, V=0.
 \endaligned \right.
 $$
The question of existence of
  suitable nontrivial solution to (ST) is equivalent to the  that of existence of
  nontrivial self-similar finite time blowing up
  solution to the transport equation. We will establish the
  following theorem.
\begin{theorem}
Suppose there exist $\a \neq -1$, $\beta \in \Bbb R$ and  solution
$(V, \Theta) $ to the system  (ST) with  $\Theta \in L^{p_1} (\Bbb
R^n )\cap L^{p_2} (\Bbb R^n)$ for some $p_1, p_2 $ such that $0<p_1
< p_2\leq \infty$.
 Then, $\Theta=0$.
\end{theorem}
This theorem  is a corollary of the following one.
\begin{theorem}
Suppose there exists $T>0$ such that there exists  a
 representation of  the solution $\theta(x,t)$ to the system (TE) by
 \bb\label{thm22}
 \t(x,t) = \Psi(t)\Theta (\Phi(t)x) \qquad \forall t\in [0, T)
 \ee
 where $\Psi(\cdot)\in C([0, T );(0, \infty))$,
 $\Phi(\cdot)\in C([0, T );\Bbb R^{n\times n})$
 with $\mathrm{det}(\Phi(t))\neq 0$ on $[0, T)$; there exists
 $p_1<p_2$ with  $p_1, p_2 \in (0, \infty]$ such that $\Theta \in L^{p_1}
 (\Bbb R^n )\cap L^{p_2} (\Bbb R^n )$.  Then,
necessarily
 either $\mathrm{det}(\Phi(t))\equiv \mathrm{det}(\Phi(0))$
  and $\Psi (t)\equiv\Psi (0)$ on $ [0, T)$, or $\Theta =0$.
\end{theorem}
{\bf Proof.} Similarly to the proof of Theorem 1.2 the
representation (\ref{thm22}) reduces to the form,
 \bb\label{transp}
  \theta (x,t) =
G(t) \theta_0 (F(t)x),
 \ee
 where $G(t)=\Psi (t)/\Psi (0)$,
$F(t)=\Phi (t)[\Phi (0)]^{-1}$. By standard $L^p$-interpolation and
the relation between $\theta_0$ and $\Theta$ by
 $\t_0 (x)=\Psi (0) \Theta (\Phi (0)x)$, we have that
 $\Theta \in L^{p_1} (\Bbb R^n )\cap L^{p_2} (\Bbb R^n )$ implies
 $\t_0 \in
L^p (\Bbb R^n )$ for all $p\in [p_1, p_2 ]$. As in the proof of
Theorem 1.2 we denote by $\{ X(a,t) \}$ and $\{ A(x,t)\}$ the
particle trajectory map and the back to label map respectively, each
one of which is defined by $v(x,t)$. As the solution of the first
equation of (TE) we have $\theta(X(a,t),t)=\theta_0 (a)$, which can
be rewritten as $\theta (x,t) = \theta_0 (A(x,t))$  in terms of the
back to label map. This, combined with (\ref{transp}),  provides us
with the relation
  \bb\label{transp1}
\theta_0 (A(x,t))=G(t) \theta_0 (F(t)x).
 \ee
Using the fact, det$(\nabla A(x,t))=1$, we compute $L^p(\Bbb R^n)$
norm of (\ref{transp1}) to have
 \bq\label{pr22}
\|\theta _0 \|_{L^p}
 &=& |G(t)| |\mathrm{det }(F(t))|^{-\frac1p} \left(\int_{\Bbb R^n}
 |\t (F(t)x)|^p |\mathrm{det}(F(t))| dx\right)^{\frac1p}\n \\
 &=&|G(t)| |\mathrm{det} (F(t))|^{-\frac1p} \|\t_0 \|_{L^p}
 \eq
 for all $t\in [0, T)$ and $p\in [p_1 , p_2]$.
 Suppose $\t_0 \neq 0$, which is equivalent to $\Theta\neq 0$,
  then we divide (\ref{pr22}) by $\|\t_0
 \|_{L^p}$ to obtain
 $|G(t)|^p =\mathrm{det} (F(t))$ for all $t\in [0, T)$ and $p\in [p_1 ,
 p_2]$, which is possible only if $G(t)=\mathrm{det} (F(t))=1$
 for all $t\in [0, T)$. Hence,
 $\Psi (t)\equiv\Psi (0),$ and det$(\Phi (t))\equiv$ det$(\Phi (0))$.  $\square$\\
\ \\
\noindent{\bf Proof of Theorem 2.1} We apply Theorem 2.2 with
$$\Phi
(t)= (T_* -t)^{-\frac{1}{\a +1}} I \quad \mbox{and}\quad \Psi (t)=
(T_* -t)^{-\beta},
$$
where $I$ is the unit matrix in $\Bbb R^{n\times n}$. Then,
$$
\mathrm{det}(\Phi (t)) =(T_* -t)^{-\frac{n}{(\a +1)}} \neq
\mathrm{det}(\Phi (0))= T_* ^{-\frac{n}{(\a+1)}} \quad
\mbox{if}\quad \a\neq -1, t\neq 0.
$$
Hence, by Theorem 2.2 we have $\Theta =0$. $\square$\\
\ \\
Below we present some examples of fluid mechanics, where we can
apply similar argument to the above to prove nonexistence of
nontrivial self-similar blowing up solutions.

\subsubsection*{A. The density-dependent Euler equations}

 The density-dependent Euler
equations in $\Bbb R^n$, $n\geq 2$, are the following system.
 \[
\mathrm{ (E_1)}
 \left\{ \aligned
 &\frac{\partial \rho v}{\partial t} +\mathrm{div}\, (\rho v\otimes v) =-\nabla
 p, \\
  &\frac{\partial \rho }{\partial t} +v\cdot \nabla \rho =0,\\
 &\quad \textrm{div }\, v =0, \\
  &v(x,0)=v_0 (x),\quad \rho (x,0)=\rho_0 (x),
  \endaligned
  \right.
  \]
  where $v=(v_1, \cdots ,v_n )=v(x,t)$ is the velocity,
  $\rho =\rho (x,t)\geq 0$ is the scalar
  density of the fluid, and $p=p(x,t)$ is the pressure.
  We refer to section 4.5 in \cite{lio} for more
detailed introduction of this system. Here we just note that this
system reduces to the homogeneous Euler system of the previous
section when $\rho\equiv 1$. The question of finite time blow-up
for the system is wide open even in the case of $n=2$, although we
have local in time existence result of the classical solution and
its finite time blow-up criterion(see e.g. \cite{bei, cha2}).
 The system $(E_1)$
has scaling property that
  if $(v, \rho, p)$ is a
solution of the system $(E_1)$, then for any $\lambda >0$ and
$\alpha \in \Bbb R $ the functions
 \bb
 \label{selfa}
  v^{\l, \a }(x,t)=\lambda ^\alpha v(\lambda x, \l^{\a +1}
  t),\quad
  \rho^{\l, \a,\beta} (x,t)= \l ^{\beta} \rho (\lambda x, \l^{\a +1}
  t),
  \ee
  \bb
  p^{\l, \a,\beta}(x,t)= \lambda ^{2\alpha +\beta }p (\lambda x,
\l^{\a +1}
  t)
  \ee
  are also solutions of $(E_1)$ with the initial data
  $$ v_0^{\l,\a}(x)=\lambda ^\alpha v_0
   (\lambda x),\quad  \rho_0^{\l,\a,\beta}(x)=
  \l ^{\beta} \rho_0(\lambda x ).
  $$
  In view of the scaling
  properties in (\ref{selfa}), we should check if there exists
  nontrivial
  solution $(v(x, t), \rho (x,t))$ of $(E_1)$ of the form,
  \bq
  \label{selfb1}
 v(x, t)&=&\frac{1}{(T_*-t)^{\frac{\a}{\a+1}}}
V\left(\frac{x}{(T_*-t)^{\frac{1}{\a+1}}}\right),\\
\label{selfb2}
 \rho  (x,t)&=&\frac{1}{(T_*-t)^{\beta}}
 R\left(\frac{x}{(T_*-t)^{\frac{1}{\a+1}}}\right)
 \eq
 for $\a \neq -1$ and $t$ sufficiently close to $T_*$. The solution $(v,\rho )$ of the form
 (\ref{selfb1})-(\ref{selfb2})
 is called the self-similar blowing up solution of the system
 $(E_1)$. The following theorem establish the
 nonexistence of nontrivial self-similar blowing up solution of the system
 $(E_1)$,  which is immediate from Theorem 2.2.
\begin{theorem}
Suppose there exist $\a \neq -1$ and solution $(v, \rho)$ to the
system $(E_1)$ of the form (\ref{selfb1})-(\ref{selfb2}), for which
there exists $p_1 , p_2$  with $0<p_1 < p_2 \leq\infty$ such that $
R\in L^{p_1} (\Bbb R^n ) \cap L^{p_2} (\Bbb R^n )$. Then, $R=0$.
\end{theorem}

\subsubsection*{ B. The 2D Boussinesq system}

The Boussinesq system for the inviscid  fluid flows in $\Bbb R^2$ is
given by
$$
  (B)\left\{
 \aligned
 &\frac{\partial v}{\partial t} +(v\cdot \nabla )v =-\nabla p + \theta e_1 ,
  \\
&\frac{\partial \theta}{\partial t} +(v\cdot \nabla )\theta =0,\\
 &\quad\textrm{div }\, v =0 ,\\
 &v(x,0)=v_0 (x), \qquad \theta(x,0)=\theta_0 (x)
 \endaligned
 \right.
 $$
where $v=(v_1, v_2)=v(x,t)$ is the velocity, $e_1 =(1,0)$, and
$p=p(x,t)$ is the pressure, while  $\theta =\theta (x,t)$ is the
temperature function. The local in time existence of solution and
the blow-up criterion of the  Beale-Kato-Majda type has been well
known(see e.g. \cite{e, cha3}). The question of finite time blow-up
is open until now. Here, we exclude the possibility of self-similar
finite time blow-up for the system.
 The system (B) has scaling property that
  if $(v, \t , p)$ is a
solution of the system (B), then for any $\lambda >0$ and $\alpha
\in \Bbb R$ the functions
 \bb
 \label{self2a}
  v^{\l, \a }(x,t)=\lambda ^\alpha v(\lambda x, \l^{\a +1}
  t),\quad
  \t ^{\l, \a } (x,t)= \l ^{2\a +1} \t (\lambda x, \l^{\a +1}
  t),
  \ee
  \bb
  p^{\l, \a }(x,t)= \lambda ^{2\alpha}p (\lambda x,
\l^{\a +1}
  t)
  \ee
  are also solutions of (B) with the initial data
  $$ v_0^{\l,\a}(x)=\lambda ^\alpha v_0
   (\lambda x),\quad  \t_0^{\l,\a}(x)=
  \l ^{2\a +1} \t_0(\lambda x ).
  $$
  In view of the scaling
  properties in (\ref{self2a}), the
 self-similar blowing-up
  solution $(v(x, t), \t (x,t))$ of (B)  should of the form,
  \bq
  \label{self2b1}
 v(x, t)&=&\frac{1}{(T_*-t)^{\frac{\a}{\a+1}}}
V\left(\frac{x}{(T_*-t)^{\frac{1}{\a+1}}}\right),\\
\label{self2b2}
 \t (x,t)&=&\frac{1}{(T_*-t)^{2\a+1}}
 \Theta \left(\frac{x}{(T_*-t)^{\frac{1}{\a+1}}}\right),
 \eq
 where $\alpha
\neq -1 $. We have the following nonexistence result of such type of
solution.
\begin{theorem}
There exists no nontrivial solution $(v, \t )$ of the system $(B)$
of the form (\ref{self2b1})-(\ref{self2b2}), if there exists $p_1,
p_2 \in (0, \infty]$, $p_1 < p_2$, such that $\Theta\in L^{p_1}
(\Bbb R^2 ) \cap L^{p_2} (\Bbb R^2 )$, and $V\in H^m (\Bbb R^2 )$,
$m>2$.
\end{theorem}
\noindent{\bf Proof.} Similarly to  the proof of Theorem 2.1, we
first conclude $\Theta =0$, and hence $\theta (\cdot ,t) \equiv 0$
on $[0,T_*)$. Then, the system (B) reduces to the 2D incompressible
Euler equations, for which we have global in time regular solution
for $v_0 \in H^m (\Bbb R^2 ), m>2$(see e.g. \cite{kat2}). Hence, we
should have
$v(\cdot ,t)\equiv 0$ on $[0,T_*)$. $\square$\\
\ \\
\noindent{\it Note added to the proof.} Similar proof to the above
leads to the nonexistence of  self-similar blowing up solution to
the axisymmetric 3D Euler equations  with swirl of the form,
(\ref{vel}), if $\Theta = rV^\theta$ satisfies the condition of
Theorem 2.4, and curl $V\in H^m (\Bbb R^3 )$, $m>5/2$, where
$r=\sqrt{x_1^2 +x_2 ^2}$, and $V^\theta$ is the angular component
of $V$. Indeed, applying Theorem 2.2 to the $\theta$-component of
the Euler equations, $\frac{D}{Dt}( rv^\theta )=0$,  we show that
$v^\theta=0$ as in the above proof, and then we use the global
regularity result for the 3D axisymmetric Euler equations without
swirl(\cite{maj0, ray}) to conclude that $(v^r , v^{3})$ is also
zero. \\

\subsubsection*{ C. The 2D quasi-geostrophic equation}

 The following 2D quasi-geostrophic
equation(QG) models the dynamics of the mixture of cold and hot air,
and the fronts between them.
$$
(QG) \left\{\aligned
& \frac{\partial \theta}{\partial t} +(v\cdot \nabla )\theta =0 , \\
&v=-\nabla ^{\bot }(-\Delta )^{-\frac{1}{2}}\theta \left( =\nabla
^{\bot }\int_{\Bbb R^2} \frac{\theta (y,t)}{|x-y|} dy\right),
\nonumber \\
&\theta (x,0)=\theta _{0}(x),
\endaligned
\right.
 $$
 where
 $\nabla ^{\bot }=(-\partial _{2},\partial _{1})$.
 Besides its physical significance, mainly due to its similar structure to
 the 3D Euler equations, there have
 been many recent studies on this system(see e.g. \cite{con3, cor1, cor2} and
 references therein).
 Although the question of
 finite time singularities is still open, some type
 of scenarios of singularities have been excluded(\cite{cor1, cor2, cor4}).
 Here we exclude the scenario of self-similar singularity.
The  system (QG) has the scaling property that
  if $\t$ is a
solution of the system, then for any $\lambda >0$ and $\alpha \in
\Bbb R $ the functions
 \bb
 \label{self3a}
  \t^{\lambda, \alpha}(x,t)=\lambda ^\alpha \t(\lambda x, \l^{\a +1}
  t)
  \ee
  are also solutions of (QG) with the initial data
  $ \t ^{\lambda, \alpha}_0(x)=\lambda ^\alpha \t_0
   (\lambda x)$.
Hence, the self-similar blowing up solution should be of the
  form,
  \bb
  \label{self3b}
 \t(x, t)=\frac{1}{(T_*-t)^{\frac{\a}{\a+1}}}
\Theta\left(\frac{x}{(T_*-t)^{\frac{1}{\a+1}}}\right)
 \ee
 for $t$ sufficiently close $T_*$ and $\a \neq -1$.
 Applying the same argument as in the proof of Theorem 2.1, we have
 the following theorem.
\begin{theorem}
There exists no nontrivial solution $ \t  $ to the system $(QG)$ of
the form (\ref{self3b}),  if there exists $p_1, p_2 \in (0, \infty
]$, $p_1 < p_2$, such that $\Theta\in L^{p_1} (\Bbb R^2 ) \cap
L^{p_2} (\Bbb R^2 )$.
\end{theorem}

\end{document}